\documentclass[12pt,oneside,reqno]{amsart}

\usepackage{amssymb}
\usepackage{txfonts}
\usepackage{bbm}
\usepackage{cases}
\usepackage{amsmath}
\usepackage{graphicx}
\usepackage{mathrsfs}
\usepackage{stmaryrd}
\usepackage{amsfonts}
\usepackage{enumerate,amsmath,amssymb,amsthm}

\numberwithin{equation}{section}

\newcommand{\be}{\begin{eqnarray}}
\newcommand{\ee}{\end{eqnarray}}
\newcommand{\ce}{\begin{eqnarray*}}
\newcommand{\de}{\end{eqnarray*}}
\newtheorem{theorem}{Theorem}[section]
\newtheorem{lemma}[theorem]{Lemma}
\newtheorem{remark}[theorem]{Remark}
\newtheorem{definition}[theorem]{Definition}
\newtheorem{proposition}[theorem]{Proposition}
\newtheorem{Examples}[theorem]{Example}
\newtheorem{corollary}[theorem]{Corollary}

\def\eps{\varepsilon}

\def\e{\mathrm{e}}

\def\p{\partial}

\def\[{{\Big[}}
\def\]{{\Big]}}
\def\<{{\langle}}
\def\>{{\rangle}}
\def\({{\Big(}}
\def\){{\Big)}}

\def\bx{{\mathbf{x}}}

\def\dif{{\mathord{{\rm d}}}}

\def\no{\nonumber}
\def\={&\!\!=\!\!&}

\def\mD{{\mathbb D}}
\def\mE{{\mathbb E}}

\def\mI{{\mathbb I}}

\def\mK{{\mathbb K}}

\def\mN{{\mathbb N}}

\def\mP{{\mathbb P}}

\def\mR{{\mathbb R}}

\def\1{{\mathbf{1}}}

\def\sL{{\mathscr L}}

\def\E{\mathbb E}

\def\geq{\geqslant}
\def\leq{\leqslant}

\def\eps{\varepsilon}

\def\e{\mathrm{e}}

\def\p{\partial}

\def\[{{\Big[}}
\def\]{{\Big]}}
\def\<{{\langle}}
\def\>{{\rangle}}
\def\({{\Big(}}
\def\){{\Big)}}

\def\bx{{\mathbf{x}}}

\def\dif{{\mathord{{\rm d}}}}

\def\no{\nonumber}
\def\={&\!\!=\!\!&}
\def\bt{\begin{theorem}}
\def\et{\end{theorem}}
\def\bl{\begin{lemma}}
\def\el{\end{lemma}}
\def\br{\begin{remark}}
\def\er{\end{remark}}
\def\bx{\begin{Examples}}
\def\ex{\end{Examples}}
\def\bd{\begin{definition}}
\def\ed{\end{definition}}
\def\bp{\begin{proposition}}
\def\ep{\end{proposition}}
\def\bc{\begin{corollary}}
\def\ec{\end{corollary}}

\def\wt{\widetilde}

\def\geq{\geqslant}
\def\leq{\leqslant}

 \def\R{\mathbb R}
 \def\R{\mathbb R}

\def\<{\langle} \def\>{\rangle}

 \def\beq{\begin{equation}}  
 
\def\e{\text{\rm{e}}}

\allowdisplaybreaks

\begin{document}

\title{Heat kernels for time-dependent non-symmetric stable-like operators}

\author{Zhen-Qing Chen and  Xicheng Zhang}

 \address{Zhen-Qing Chen:
Department of Mathematics, University of Washington, Seattle, WA 98195, USA\\
Email: zqchen@uw.edu
 }

\address{Xicheng Zhang:
School of Mathematics and Statistics, Wuhan University,
Wuhan, Hubei 430072, P.R.China\\
Email: XichengZhang@gmail.com
 }

\thanks{
 This work is partially supported by  Simons Foundation grant 520542  and by an NNSFC grant of China (No. 11731009). }

\begin{abstract}
When studying non-symmetric nonlocal operators 
$$
\sL f(x) =  \int_{\R^d}
\left( f(x+z)-f(x)-\nabla f(x)\cdot z 1_{\{|z|\leq 1\}} \right) \frac{\kappa (x, z)}{|z|^{d+\alpha}} \dif z , 
$$
where $0<\alpha<2$ and $\kappa (x, z)$ is a function on $\R^d\times \R^d$ that is bounded between
two positive constants, it is customary to assume that $\kappa (x, z)$ is symmetric in $z$. 
In this paper,  we study heat kernel of $\sL$ and derive its two-sided sharp bounds without
 the  symmetric assumption $\kappa(x,z)=\kappa(x,-z)$.
 In fact,  we allow the kernel $\kappa$ to be time-dependent   and also derive    gradient estimate  
 when $\beta\in(0\vee (1-\alpha),1)$ as well as  fractional derivative estimate of order  $\theta\in(0,(\alpha+\beta)\wedge 2)$ for the heat kernel, 
 where $\beta$ is 
the H\"older index of $x\mapsto\kappa(x,z)$. Moreover, when $\alpha\in(1,2)$, the drift perturbation with drift in  Kato's class 
is also considered.
As an application, when $\kappa(x,z)=\kappa(z)$ does not depend on $x$, we show the boundedness of nonlocal Riesz's transorfmation:
for any $p>2d/(d+2\alpha)$,
$$
\|\sL^{1/2}f\|_p\asymp \|\Gamma(f)^{1/2}\|_p,
$$
where $\Gamma(f):=\frac{1}{2}\sL (f^2)-f\sL f$ is the carr\'e du champ operator associated with $\sL$,  
and $\sL^{1/2}$ is the square root operator of $\sL$ defined by using Bochner's subordination. Here 
$\asymp$ means that both sides are comparable up to a 
 constant multiple.

\bigskip

	
	\medskip


\end{abstract}

\maketitle

\section{Introduction}

For $\alpha\in(0,2)$, we consider the following nonlocal and non-symmetric operator:
\begin{align}\label{Non}
\sL^\kappa_t f(x):=\int_{\mR^d} \left(f(x+z)-f(x)-z^{(\alpha)}\cdot\nabla f(x)\right)
\frac{\kappa(t,x,z)}{|z|^{d+\alpha}}\dif z,
\end{align}
where 
\begin{equation}\label{e:1.2}
z^{(\alpha)}:= \left(1_{\alpha\in(1,2)}+1_{|z|\leq 1}1_{\alpha=1} \right) z,
\end{equation}
and
$\kappa(t,x,z):\mR_+\times\mR^d\times\mR^d\to\mR$ satisfies that for some $\kappa_0>1$, $\beta\in(0,1]$ and $\beta'\geq 0$,
\begin{align}\label{Con1}
\kappa^{-1}_0\leq \kappa(t,x,z)\leq \kappa_0,\ \ |\kappa(t,x,z)-\kappa(t,y,z)|\leq \kappa_0|x-y|^\beta(1+|z|^{\beta'}),
\end{align}
and for $\alpha=1$,
\begin{align}\label{Con2}
\int_{R_0\leq |z|\leq R_1}\frac{ z \kappa(t,x,z)}{|z|^{d+1}}\dif z=0
\quad \hbox{for any }  0<R_0<R_1<\infty.
\end{align}
Note that this condition is equivalent to
\begin{equation}\label{e:1.3}
\int_{  |z|\leq r} z \kappa(t,x,z)  \dif z=0  \quad \hbox{for every } r>0.
\end{equation}
The reason we use $z^{(\alpha)}$ instead of the more common $z 1_{\{|z| \leq 1\}}$
in the first order correction term in \eqref{Non} together with condition \eqref{Con2} 
is that 
this is the form 
for general $\alpha$-stable L\'evy processes on $\R^d$ when $\kappa (t, x, z)$ is a constant.

When $\kappa(t,x,z)=\kappa(x,z)$ is time-independent and $\kappa(x,z)=\kappa(x,-z)$, the heat kernel of $\sL^\kappa_t$ has been constructed in \cite{Ch-Zh}.
Moreover, sharp two-sided estimates,  gradient estimate and fractional derivative estimate of the heat kernel are obtained \cite{Ch-Zh}. 
The main goal of this paper is to  drop the symmetric assumption in $z$ and extend it to time-dependent case at the same time. It in particular gives sharp two-sided heat kernel estimates to non-symmetric $\alpha$-stable processes whose L\'evy measure is comparable to that of isotropic $\alpha$-stable process on $\R^d$.
The study of heat kernels and their estimates is an active research area in analysis and in probability theory. 
We refer the reader to the Introduction   of \cite{Ch-Zh} for a brief history on the study of heat kernels for nonlocal operators.

To state our main results, we introduce the following notations for later use: for $\beta\geq 0$ and $\gamma\in\mR$,
$$
\varrho^\beta_\gamma(t,x):=\frac{t^{\gamma/\alpha}(|x|^\beta\wedge 1)}
{(t^{1/\alpha}+|x|)^{d+\alpha}},
$$
and for any $T\in(0,\infty]$ and $\eps\in[0,T)$,
$$
\mD^T_{\eps}:=\Big\{(t,x;s,y): x,y\in\mR^d \mbox{ and } s,t\geq 0 \mbox{ with } \eps<s-t<T\Big\}.
$$
It is well known that the transition density function $p(t, x)$
for an isotropic $\alpha$-stable process on $\R^d$ has the property 
$$
p (t, x)\asymp   \varrho^0_\alpha( t,x)\quad \hbox{for all } t>0 \hbox{ and } x\in \R^d.
$$
Here and below $\asymp$ means that both sides are comparable up to  a constant multiple.

In this paper, we show the following.

\bt\label{heat}
Under \eqref{Con1} and \eqref{Con2}, there is a unique continuous function $p^\kappa_{t,s}(x,y)$ on $\mD^\infty_0$
(called the fundamental solution or heat kernel of $\sL^\kappa_t$) satisfying
\begin{align}\label{EQ}
\p_t p^\kappa_{t,s}(x,y)+\sL^\kappa_t p^\kappa_{t,s}(\cdot,y)(x)=0
\end{align}
for Lebesgue almost all $t\in[0,s]$, and 
\begin{enumerate}[\rm (i)]
\item (Uniform continuity)
For any bounded and uniformly continuous function $f(x)$, 
      \begin{align}
        \lim_{|t-s|\to 0}\|P^\kappa_{t,s}f-f\|_\infty=0,\label{cz1}
      \end{align}
where $P^\kappa_{t,s}f(x):=\int_{\mR^d}p^\kappa_{t,s}(x,y)f(y)\dif y$.

\item There exits some $\theta >\alpha $ so that 
$(t,x)\mapsto \Delta^{\theta/2}  P^\kappa_{t, s}f (x)$   is bounded and continuous on $[0,s-\eps]\times\mR^d$ for every $s>\eps>0$.
\end{enumerate}
Moreover, $p^\kappa_{t,s}(x,y)$ enjoys the following properties:
\begin{enumerate}[(i)]
\item[\rm (iii)] (Two-sided estimate) For any $T>0$, there is a constant 
$c_0=c_0(d, \alpha, \beta, \beta', \kappa_0)\geq 1$ 
 such that on $\mD^T_0$,
\begin{align}\label{GR1}
c^{-1}_0\varrho^0_\alpha(s-t,x-y)\leq p^{\kappa}_{t,s}(x,y)\leq c_0 \varrho^0_\alpha(s-t,x-y).
\end{align}

\item[\rm (iv)](Fractional derivative estimate) 
There exits an $\eps=\eps(\beta,\beta',\alpha)\in(0,2-\alpha)$ such that for any $\theta\in[0,\alpha+\eps)$,
$(t,x)\mapsto \Delta^{\theta/2}_xp^\kappa_{t,s}(\cdot,y)(x)$ is continuous on $[0,s)\times\mR^d$, and for any $T>0$, there is a constant $c_1>0$ such that   on $\mD^T_0$,
\begin{align}\label{Grad}
|\Delta^{\theta/2}_xp^\kappa_{t,s}(\cdot,y)(x)|\leq c_1\varrho^0_{\alpha-\theta}(s-t,x-y).
\end{align}

\item[\rm (v)]
(Gradient estimate) If $\alpha\in[1,2)$ or $\alpha+\beta>1$ and $\beta'=0$ in \eqref{Con1}, 
      then for any $T>0$, there is a constant $c_2>0$ such that  on $\mD^T_0$,
\begin{align}\label{GR2}
|\nabla_x \log p^\kappa_{t,s}(\cdot,y)(x)|\leq c_2 (s-t)^{-1/\alpha}.
\end{align}

\item[\rm (vi)]
(Conservativeness)
For every $0<s<t$ and $x\in \mR^d$, 
\begin{align}
\int_{\mR^d}p^\kappa_{t,s} (x,y)\dif y=1.\label{EK1}
\end{align}

\item[\rm (vii)] (C-K equation) 
For all $0<t<r<s<\infty$ and $x,y\in\mR^d$,
the following Chapman-Kolmogorov equation holds:
\begin{align}\label{eq21}
\int_{\mR^d}p^\kappa_{t,r} (x,z)p^\kappa_{r,s} (z,y)\dif z=p^\kappa_{t,s} (x,y).
\end{align}

\item[\rm (viii)] (Generator) 
For any $f\in C_b^2(\mR^d)$, we have
      \begin{align}
        P^\kappa_{t,s}f(x)-f(x)=\int^s_t\!P^\kappa_{t,r}\sL^\kappa_rf(x)\dif r=\int^s_t\!\sL^\kappa_r P^\kappa_{r,s}f(x)\dif r.\label{eqge}
      \end{align}
\end{enumerate}
\et

\br  \label{R:1.2} \rm \begin{enumerate}[\rm (i)] 
\item  It should be noticed that in \eqref{Grad}, the differentiability index $\theta$ could be greater than $\alpha$.
From the proof below, one sees that if $\beta'=0$ in \eqref{Con1}, then we can take $\eps=(2-\alpha)\wedge\beta$ in 
 {\rm (iv)} of Theorem \ref{heat}. 
 
 \item It follows from Theorem \ref{heat} that the fundamental solution $p^\kappa_{t, s}(x, y)$ uniquely determines
 a time-inhomogeneous Feller process $X:=\{X^{t, x}_s; s\geq 0; \mP^{s, x},  s\geq 0, x\in \mR^d\}$ that has $p^\kappa_{t, s}(x, y)$
 as its transition density function with respect to the Lebesgue measure on $\mR^d$. Clearly, by Theorem \ref{heat}(viii), 
 the probability law of $X$ solves the martingale
 problem for $(\sL^\kappa_r, C^\infty_c  (\mR^d))$ in the sense that for every $f\in C^\infty_c (\mR^d)$ and $(t, x)\in [0, \infty) \times \mR^d$, 
 $$
 M^f_s :=f(X^{t,x}_s)- f(x) -\int_t^s \sL^\kappa_r f (X^{t,x}_r) \dif r; \quad s\geq t, 
 $$
 is an $\mP^{t,x}$-martingale. Here $C^\infty_c(\mR^d))$ is the space of smooth functions on $\mR^d$ with compact support.
 This gives a constructive proof of the existence of solution to the martingale problem $(\sL^\kappa_r, C^\infty_c  (\mR^d))$. 
 The existence and uniqueness of solutions to  the martingale problem  for  $(\sL^\kappa_r, C^\infty_c  (\mR^d))$ have been established in 
 \cite[Proposition 3]{MP} and \cite[Theorem 4.6]{Ch-Zh0}.  
 \end{enumerate}
\er

\medskip

Next we consider the perturbation of $\sL^\kappa_t$ by a drift $b$ belonging to some Kato's class when $\alpha\in(1,2)$. First of all, we introduce the following 
Kato's class as in \cite{Wa-Zh}.
\bd
For $\alpha\in(1,2)$, a Borel measurable function $f: [0,\infty)\times\mR^d\to\mR$ is 
said to be 
 in Kato's class $\mK^\alpha_d$ if
$$
\lim_{\eps\downarrow 0}K^\alpha_f(\eps)=0,
$$
where
$$
K^\alpha_f(\eps):=\eps  \sup_{(t,x)\in [0,\infty)\times\mR^d} \int^\eps_0\!\!\!\int_{\mR^d}  
 \frac{\varrho^0_\alpha(s,x-y)|f(t\pm s, y)|}{s^{1/\alpha}
 (\eps-s)^{1/\alpha}}\dif y\dif s.
$$
Here we have extended $f$ to $\mR$ by setting $f(t,\cdot)=0$ for $t<0$.
\ed
\br
By H\"older's inequality, one sees that $L^q_{loc}
(\mR_+; L^p (\mR^d))\subset\mK^\alpha_d$ 
provided that $\frac{d}{p}+\frac{\alpha}{q}<\alpha-1$.
In particular, $\mK^\alpha_d$ contains all bounded functions on $\mR_+\times \mR^d$.
\er

When $\alpha \in (1, 2)$ and $b(t, x)\in \mK^\alpha_d$, $\sL^\kappa_t + b(t, x)\cdot \nabla$
can be viewed as the perturbation of $\sL^\kappa_t$ by $b(t, x)\cdot \nabla$.
Hence, heuristically, the fundamental solution $p^{\kappa, b}_{t, s}(x, y)$ of 
$\sL^\kappa_t + b(t, x)\cdot \nabla$ is related to $p^\kappa_{t, s}(x, y)$ of $\sL^\kappa_t$
by the following Duhamel's formula: for any $t<s$ and $x, y\in \mR^d$, 
\begin{equation}\label{e:1.13}
p^{\kappa, b}_{t, s}(x, y) = p^{\kappa }_{t, s}(x, y) + \int_t^s\!\!\! \int_{\mR^d}
p^{\kappa, b}_{t, r}(x, z) b(r, z) \nabla_z p^{\kappa }_{r, s}(z, y) \dif z \dif r. 
\end{equation}

The following result can be shown as in \cite{Bo-Ja} and 
\cite{Wa-Zh}, while the uniqueness can be shown as in \cite[Theorem 3.10]{CW} or
\cite[Theorem 1.1]{CHXZ}. We omit the details.

\bt \label{T:1.5}
Let $\alpha\in(1,2)$ and $b\in\mK^\alpha_d$. Under \eqref{Con1} and \eqref{Con2}, 
there is a unique continuous function $p^{\kappa,b}_{t,s}(x,y)$ on $\mD^\infty_0$ 
satisfying \eqref{e:1.13} and that for every $T>0$, there is a constant $c_0>1$ such that on $\mD^T_0$,
\begin{align*} 
  |p^{\kappa,b}_{t,s}(x,y)|
  \leq c_0 \, \varrho^0_\alpha(s-t,x-y).
\end{align*}
Moreover, $p^{\kappa,b}_{t,s}(x,y)$ enjoys the following properties.

\begin{enumerate}[\rm (i)]
\item (Two-sided estimate) For any $T>0$, there is a constant $c_1>1$ such that on $\mD^T_0$,
\begin{align*}
c^{-1}_1\varrho^0_\alpha(s-t,x-y)\leq p^{\kappa,b}_{t,s}(x,y)\leq c_1 \varrho^0_\alpha(s-t,x-y).
\end{align*}
\item (Conservativeness)
For every $0<s<t$ and $x\in \mR^d$, 
\begin{align*}
\int_{\mR^d}p^{\kappa,b}_{t,s} (x,y)\dif y=1.
\end{align*}
\item (C-K equation) For all $0<t<r<s<\infty$ and $x,y\in\mR^d$,
the following Chapman-Kolmogorov equation holds:
\begin{align*}
\int_{\mR^d}p^{\kappa,b}_{t,r} (x,z)p^{\kappa,b}_{r,s} (z,y)\dif z=p^{\kappa,b}_{t,s} (x,y).
\end{align*}
 \item (Generator) For any $f\in C_b^2(\mR^d)$, we have
      \begin{align*}
           P^{\kappa,b}_{t,s}f(x)-f(x)=\int^s_t\!P^{\kappa,b}_{t,r}\sL^{\kappa,b}_rf(x)\dif r=\int^s_t\!\sL^{\kappa,b}_rP^{\kappa,b}_{r,s}f(x)\dif r,
      \end{align*}
    where  $P^{\kappa,b}_{t,s}f(x):=\int_{\mR^d}p^{\kappa,b}_{t,s}(x,y)f(y)\dif y$ and 
       $\sL^{\kappa,b}_r:=\sL^\kappa_r+b(r,\cdot)\cdot\nabla$.
  \item (Uniformly continuity) For any bounded and uniformly continuous function $f(x)$,
      \begin{align*}
        \lim_{|t-s|\to 0}\|P^{\kappa,b}_{t,s}f-f\|_\infty=0.
      \end{align*}
      \item (Gradient estimate) For any $T>0$, there is a constant $c_2>0$ such that  on $\mD^T_0$,
\begin{align*}
|\nabla_x \log p^{\kappa,b}_{t,s}(\cdot,y)(x)|\leq c_2 (s-t)^{-1/\alpha}.
\end{align*}
\end{enumerate}
\et

\medskip

\br \rm Under the condition of Theorem \ref{T:1.5},   the fundamental solution $p^{\kappa, b}_{t, s}(x, y)$ uniquely determines
 a time-inhomogeneous Feller process $X:=\{X^{t, x}_s; s\geq 0; \mP^{s, x},  s\geq 0, x\in \mR^d\}$ that has $p^\kappa_{t, s}(x, y)$
 as its transition density function with respect to the Lebesgue measure on $\mR^d$.  It follows from Theorem \ref{T:1.5}(iv) that 
 the probability law of $X$ solves the martingale
 problem for $(\sL^{\kappa, b}_r, C^\infty_c  (\mR^d))$. When $b$ is bounded, 
 the existence and uniqueness of solutions to  the martingale problem  for  $(\sL^{\kappa, b}_r, C^\infty_c  (\mR^d))$ have been established in 
 \cite[Proposition 3]{MP} and \cite[Theorem 4.6]{Ch-Zh0}. 
 Using the approach in \cite{CLW}, one can in fact establish the uniqueness of the martingale problem for 
 $(\sL^{\kappa, b}_r, C^\infty_c  (\mR^d))$ for general $b\in\mK^\alpha_d$.  
\er

\medskip

Finally, as an application of heat kernel estimate obtained in Theorem \ref{Th1} below, we have the following boundedness of nonlocal Riesz's 
transform. 
We use $\sL^{1/2}$ to denote the ``square root" operator of $\sL$ through Bochner's subordination as follows. 
 Suppose $\kappa(t,x,z)=\kappa(z)$. Then $\sL$  defined by \eqref{Non} is the generator of a   (possibly non-symmetric) L\'evy process $X_t$
 with transition semigroup $\{P_t; t\geq 0\}$. Let $S_t$ be an independent $\frac12$-subordinator. Clearly $Y_t:=X_{S_t}$ is again a L\'evy process,
 whose generator we denote as $\sL^{1/2}$. It is well-known  (see \cite[p.216, (32.11)]{Sa} that 
 \begin{align}\label{UR1}
 \sL^{1/2}f(x) =\frac{1}{2\Gamma(1/2)}\int^\infty_0 (P_tf(x)-f(x))t^{-3/2}\dif t.
 \end{align}
 
\bt\label{Th01}
Let $\alpha\in(0,2)$ and $\sL$ be defined as in \eqref{Non} with $\kappa(t,x,z)=\kappa(z)$. Suppose
\begin{align}\label{Con3}
\kappa^{-1}_0\leq\kappa(z)\leq\kappa_0
\quad \hbox{and} \quad 
 1_{\alpha=1}\int_{|z|\leq r}z\kappa(z)\dif z=0 \hbox{ for every } r>0.
\end{align}
Then for any $f\in C^\infty_c(\mR^d)$ and $p>2d/(d+2\alpha)$,
\begin{align}\label{Riesz}
\|\sL^{1/2} f\|_p\asymp\|\Gamma(f)^{1/2}\|_p,
\end{align}
where 
 $\Gamma(f):=\frac{1}{2}\sL (f^2)-f \sL f =\int_{\mR^d}(f(\cdot+z)-f(\cdot))^2\kappa(z) /|z|^{d+\alpha} \, \dif z$.
\et

Classical Riesz's transform 
says that for any $p>1$, there is a constant $c>0$ such that (see \cite{St})
\begin{align}\label{EU4}
\|\Delta^{1/2}f\|_p\asymp\|\nabla f\|_p.
\end{align}
Up to now, there are a large amount of literatures devoting to the study of various Riesz's transformation.
Here we only 
only mention that P.A. Meyer \cite{Me} used the probabilistic technique to prove \eqref{EU4}, 
 and consider the more general problem
(called Meyer's problem by now): determining whether for nice $f$,
\begin{align}\label{EU5}
\|A^{1/2}f\|_p\asymp\|\Gamma(f)^{1/2}\|_p,\ \ p\in(1,\infty),
\end{align}
where $A$ is an abstract symmetric Markov operator and $\Gamma(f):=\frac{1}{2}A(f^2)-fAf$ is the associated carr\'e du champ operator (if it exists).
Meyer \cite{Me1} 
established \eqref{EU5}
for Ornstein-Uhlenbeck operator on Wiener space. 
Bakry \cite{Ba} later showed 
it holds for diffusion operators 
on Riemannian manifold under the condition that Ricci curvature  is bounded from below. 
It seems to us that Theorem \ref{Th01}
is the first result on Riesz's transform  
for {\it non-symmetric and nonlocal} operators.
 An interesting  open problem is  whether \eqref{Riesz} holds for nonlocal operator $\sL$ with  spatial dependent kernel $\kappa(x,z)$. 
This is a quite challenging problem even for second order elliptic operators with variable coefficients,  see, for example,
\cite{Sh} and the references therein.  We plan to  investigate  this problem in a future project.

\medskip

The main part of this paper was reported at the IMS-China International Conference on Statistics and Probability
held from June 28-July 2, 2017 at Nanning, China.
At the time when we are finalizing this paper, we notice a preprint \cite{J} by Peng Jin just posted on arXiv, 
where heat kernels for time-independent $\sL$ is studied using the approach from our previous work \cite{Ch-Zh}.
The results in \cite{J} overlap some of ours in Theorems \ref{heat} and \ref{T:1.5} in the time-independent and bounded drift  case.

\section{Heat kernel estimates of $\sL^\kappa_t$ with $\kappa(t,x,z)=\kappa(t,z)$}

Throughout this section, we assume
\begin{equation}\label{Ass1a}
\kappa^{-1}_0\leq \kappa(t,z)\leq \kappa_0
\end{equation}
and when $\alpha =1$, 
\begin{equation}\label{Ass1b}
\int_{R_0\leq |z|\leq R_1}\frac{z \, \kappa(t,z)}{|z|^{d+1}}\dif z=0 \quad \hbox{for every }
 0<R_0<R_1<\infty . 
\end{equation}

Let $N(\dif t,\dif z)$ be a time-inhomogenous Poisson random measure  
with intensity measure $\frac{\kappa(t,z)}{|z|^{d+\alpha}}\dif z\dif t$.
Define
\begin{equation}\label{e:2.3a}
X^\kappa_{t,s}:=\int^s_t\!\!\!\int_{\mR^d}z \tilde N(\dif r,\dif z)+\int^s_t\!\!\!\int_{\mR^d}(z-z^{(\alpha)})\tfrac{\kappa(r,z)}{|z|^{d+\alpha}}\dif z\dif r,
\end{equation}
where $\tilde N(\dif t,\dif z):=N(\dif t,\dif z)-\frac{\kappa(t,z)}{|z|^{d+\alpha}}\dif z\dif t$.
By It\^o's formula, we have
$$
\mE f(X^\kappa_{t,s})=\mE\int^s_t\sL^{\kappa}_r f(X^\kappa_{t,r})\dif r,\ \ f\in C^2_b(\mR^d).
$$
In particular, if we take $f(x)=\e^{\mathrm{i}\xi\cdot x}$, then one finds that  the characteristic function of $X^\kappa_{t,s}$ is given by
$$
\mE \e^{\mathrm{i}\xi\cdot X^\kappa_{t,s}}=\exp\left\{\int^s_t\!\!\!\int_{\mR^d}(\e^{\mathrm{i}\xi\cdot z}-1-\mathrm{i}\xi\cdot z^{(\alpha)})
\tfrac{\kappa(r,z)}{|z|^{d+\alpha}}\dif z\dif r\right\}.
$$
By the definition of $z^{(\alpha)}$ and a change of variable, 
we conclude from the last display that for every $\lambda >0$,
\begin{equation}\label{e:2.4}
\left\{ \lambda^{-1/\alpha} X^\kappa_{\lambda t, \lambda s} , s>t \right\} \ 
\hbox{ has the same distribution as } \
\left\{ X^{\wt \kappa}_{t, s}, s>t \right\} 
\end{equation}
when $\alpha \not=1$, where $\wt \kappa (r, z)= \kappa (\lambda r, \lambda^{1/\alpha} z)$.
This is the reason why we define $z^{(\alpha)}$ in this way as \eqref{e:1.2}. See Remark \ref{R:2.7} below.
The scaling property \eqref{e:2.4} holds for $\alpha =1$ as well but under condition 
\eqref{Ass1b}.
By the change of variables, we can write
\begin{align}\label{Ch}
\mE \e^{\mathrm{i}\xi\cdot X^\kappa_{t,s}}
=\exp\left((s-t)\int^1_0\!\!\!\int_{\mR^d}(\e^{\mathrm{i}\xi\cdot z}-1-\mathrm{i}\xi\cdot z^{(\alpha)})\tfrac{\kappa(t+(s-t)r,z)}{|z|^{d+\alpha}}\dif z\dif r\right).
\end{align}
By \eqref{Ass1a},  there is a constant  $c>0$   depending only 
on $\kappa_0,d,\alpha$ so that 
\begin{eqnarray*}
|\mE \e^{\mathrm{i}\xi\cdot X^\kappa_{t,s}}|
\leq \exp\left((s-t)\int^1_0\!\!\!\int_{\mR^d} \left(\cos (\xi\cdot z) -1 \right)
 \tfrac{\kappa(t+(s-t)r,z)}{|z|^{d+\alpha}}\dif z\dif r\right)\leq \e^{-c(t-s)|\xi|^\alpha}.
\end{eqnarray*}
Hence, $X^\kappa_{t,s}$ admits a continuous density $p^\kappa_{t,s}(x)$ given by Fourier's inverse transform
\begin{align}\label{Den1}
p^\kappa_{t,s}(x)=\int_{\mR^d}\e^{-\mathrm{i}x\cdot\xi}\mE \e^{\mathrm{i}\xi\cdot X^\kappa_{t,s}}\dif\xi
=\int_{\mR^d}\mE \e^{\mathrm{i}\xi\cdot (X^\kappa_{t,s}-x)}\dif\xi.
\end{align}
Moreover, we also have
$$
\p_tp^\kappa_{t,s}(x)+\sL^\kappa_tp^\kappa_{t,s}(x)=0 \ \hbox{ for }  s>t 
\  \hbox{ with } \ \lim_{s\downarrow t}p^\kappa_{t,s}(x) \dif x =\delta_0 (\dif x),
$$
where the limit is taken in the weak sense.

The following is the main result of this section. 

\bt\label{Th1}
Under \eqref{Ass1a} and \eqref{Ass1b}, there is a constant $c_0>1$ only depending on $\kappa_0,d,\alpha$ such that for all $t<s$ and $x\in\mR^d$,
\begin{align}\label{ER11}
c^{-1}_0\varrho^0_\alpha(s-t,x)\leq p^\kappa_{t,s}(x)\leq c_0\varrho^0_\alpha(s-t,x).
\end{align}
\et

\medskip

Notice that 
by \eqref{Ass1b} and \eqref{e:2.4}, 
\begin{align}\label{Sca}
p^\kappa_{t,s}(x)=(s-t)^{-d/\alpha}p^{\tilde\kappa}_{0,1}((s-t)^{-1/\alpha}x),
\end{align}
where $\tilde\kappa(r,z):=\kappa(s+(t-s)r,(t-s)^{1/\alpha}z)$.
Thus, to prove \eqref{ER11}, it suffices to show it for $t=0$ and $s=1$.
The main point for us is to show that $c_0$ only depends on $\kappa_0,d,\alpha$.
Let $\chi_1$ and $\chi_2$ 
be the small and large jump parts of $X^\kappa_{0,1}$ defined respectively by
\begin{align*}
&\chi_1
:=\int^1_0\!\!\!\int_{|z|\leq 1}z \tilde N(\dif r,\dif z)+\int^1_0\!\!\!\int_{|z|\leq 1}(z-z^{(\alpha)})\tfrac{\kappa(r,z)}{|z|^{d+\alpha}}\dif z\dif r,\\
&\chi_2 
:=\int^1_0\!\!\!\int_{|z|>1}z \tilde N(\dif r,\dif z)+\int^1_0\!\!\!\int_{|z|>1}(z-z^{(\alpha)})\tfrac{\kappa(r,z)}{|z|^{d+\alpha}}\dif z\dif r.
\end{align*}
Note that $\chi_1$ and $\chi_2$
are independent and have the characteristic functions
\begin{align}
\mE \e^{\mathrm{i}\xi\cdot \chi_1}&=\exp\left\{\int^1_0\!\!\!\int_{|z|\leq 1}(\e^{\mathrm{i}\xi\cdot z}-1-\mathrm{i}\xi\cdot z^{(\alpha)})
\tfrac{\kappa(r,z)}{|z|^{d+\alpha}}\dif z\dif r\right\}=:\e^{\phi_1(\xi)},\label{Den}\\
\mE \e^{\mathrm{i}\xi\cdot \chi_2}&=\exp\left\{\int^1_0\!\!\!\int_{|z|>1}(\e^{\mathrm{i}\xi\cdot z}-1-\mathrm{i}\xi\cdot z^{(\alpha)})
\tfrac{\kappa(r,z)}{|z|^{d+\alpha}}\dif z\dif r\right\}=:\e^{\phi_2(\xi)}.\label{Den0}
\end{align}
In particular, $\chi_1$ has a smooth density $p_1 (x)$ with respect to the Lebesgue measure on
$\R^d$. Since $X^\kappa_{0, 1}$ is the independent sum of $\chi_1$ and $\chi_2$, we have
\begin{align}\label{ER5}
p^\kappa_{0,1}(x)=\mE \left[ p_1 (x-\chi_2) \right].
\end{align}
To show the two-sided estimate of $p^\kappa_{0,1}(x)$, we prepare the following two lemmas.

\bl\label{Le22}
{\rm (i)} For any $R>0$, there is a $\delta>0$ only depending on $R,\kappa_0,d,\alpha$ such that
$$
\inf_{x\in B_R}p_1 (x)\geq\delta.
$$

{\rm (ii)} For any integer $m\geq 1$, there is a constant $c=c(\kappa_0,d,\alpha,m)>0$ such that
$$
p_1(x)\leq c(1+|x|)^{-m},\ x\in\mR^d.
$$
\el

\begin{proof}
(i) Let $\chi_{11}$ and $\chi_{12}$ be two independent random variables with the characteristic functions
\begin{align}
\mE \e^{\mathrm{i}\xi\cdot \chi_{11}}&=\exp\left\{\int_{|z|\leq 1}(\e^{\mathrm{i}\xi\cdot z}-1-\mathrm{i}\xi\cdot z^{(\alpha)})
\tfrac{\kappa(z)-\kappa_0/2}{|z|^{d+\alpha}}\dif z\right\}=:\e^{\phi_{11}(\xi)},\label{ER1}\\
\mE \e^{\mathrm{i}\xi\cdot \chi_{12}}&=\exp\left\{\int_{|z|\leq 1}(\e^{\mathrm{i}\xi\cdot z}-1-\mathrm{i}\xi\cdot z^{(\alpha)})
\tfrac{\kappa_0/2}{|z|^{d+\alpha}}\dif z\right\}=:\e^{\phi_{12}(\xi)},\no
\end{align}
where $\kappa(z):=\int^1_0\kappa(r,z)\dif r\geq\kappa_0$ by \eqref{Ass1a}.
Let $p_{11}$ and $p_{12}$ be the continuous distribution density of $\chi_{11}$ and $\chi_{12}$.
 Clearly, we have
\begin{align}\label{ER2}
p_1(x)=\int_{\mR^d}p_{11}(x-z)p_{12}(z)\dif z.
\end{align}
Since $\chi_{12}$ is a truncated $\alpha$-stable random variable, it is well known that $p_{12}$ is strictly positive.
On the other hand, by \eqref{ER1} we also have
$$
\mE|\chi_{11}|\leq c(\kappa_0,p,\alpha).
$$
Hence, by \eqref{ER2}, we have for any $R'>0$,
\begin{align*}
p_1(x)&=\int_{\mR^d}p_{11}(x-z)p_{12}(z)\dif z\geq\inf_{z\in B_{R'}}p_{12}(z)\int_{|z|\leq R'}p_{11}(x-z)\dif z\\
&=\inf_{z\in B_{R'}}p_{12}(z)\Big(1-\mP(|\chi_{11}-x|>R')\Big)\\
&\geq\inf_{z\in B_{R'}}p_{12}(z)\Big(1-\mE|\chi_{11}-x|/R'\Big),
\end{align*}
which yields (i) by choosing $R'$ large enough.

(ii) Using Fourier's inverse  tansform, 
for every integer $m\geq 1$,  we have by \eqref{Den}
that
\begin{eqnarray*}
(1+|x|^2)^m p_1 (x)
\leq  (2\pi)^d \int_{\R^d}  \left|  (\mI-\Delta)^m \e^{\phi_1(\xi)} \right| \dif \xi  <\infty.
\end{eqnarray*}
The proof is complete.
\end{proof}
 
\bl\label{Le23}
For any $R>2$, there is a constant $c_1=c_1(R,\kappa_0,d,\alpha)>0$ such that
\begin{align}\label{EQ1}
c_1^{-1}(1+|x|)^{-d-\alpha}\leq \mP(\chi_2\in B_R(x))\leq c_1(1+|x|)^{-d-\alpha}.
\end{align}
\el

\begin{proof}
Observe that by \eqref{Den0},
\begin{align*}
\mE \e^{\mathrm{i}\xi\cdot \chi_2}&=\exp\left(\int_{\mR^d}(\e^{\mathrm{i}\xi\cdot z}-1)\nu(\dif z)\right)
\exp\left( -\mathrm{i}\xi\cdot b\right),
\end{align*}
where $\nu(\dif z):=1_{\{|z|>1\}}|z|^{-d-\alpha}\left(\int^1_0\kappa(r,z)\dif r\right)\dif z$ and $b:=\int_{\mR^d}z^{(\alpha)}\nu(\dif z)$.
Let $\eta:=\{\eta_n,n\in\mN\}$ be a family of i.i.d. random variables in $\mR^d$ with distribution $\nu/\lambda$,  
where 
$$
\lambda:=\nu(\mR^d)\leq \kappa_0 \int_{|z|>1} |z|^{-d-\alpha} \dif z<\infty.
$$ 
Let $S_0=0$ and 
$S_n:=\eta_1+\cdots+\eta_n.$
Let $N$ be a Poisson random variable with parameter $\lambda$, which is independent of $\eta$.
It is easy to see that
$$
S_N\stackrel{(d)}{=}\chi_2+b.
$$
Now, by definition we have
\begin{align*}
\mP(\chi_2\in B_R(x))&=\mP(S_N\in B_R(x+b))
=\sum_{n=1}^\infty\mP\Big(S_n\in B_R(x+b)\Big)\mP(N=n)\\
&=\e^{-\lambda}\sum_{n=1}^\infty\frac{1}{n!}\int_{\mR^{nd}}1_{\sum_{j=1}^n z_j\in B_R(x+b)}\nu(\dif z_1)\cdots\nu(\dif z_n).
\end{align*} 
When $|x+b|<R+1$, 
the upper bound in \eqref{EQ1} for 
$ \mP(\chi_2\in B_R(x))$ trivially holds. 
Thus we assume that $|x+b|\geq R+1$.
Notice that $\sum_{j=1}^n z_j\in B_R(x+b)$ implies that there is at least one $i$ such that $|z_i|>(|x+b|-R)/n$. Hence,
$$
\mP(\chi_2\in B_R(x))\leq\e^{-\lambda}\sum_{n=1}^\infty\frac{1}{n!}
\left(\sum_{i=1}^n \int_{\mR^{nd}}1_{\sum_{j=1}^n z_j\in B_R(x+b)}1_{|z_i|>(|x+b|-R)/n}\nu(\dif z_1)\cdots\nu(\dif z_n)\right).
$$
Recalling $\nu(\dif z_i)=1_{|z_i|>1}|z_i|^{-d-\alpha}\left(\int^1_0\kappa(r,z_i)\dif r\right)\dif z_i$ and by \eqref{Ass1a}, we get
\begin{eqnarray*}
&& \mP(\chi_2\in B_R(x)) \\
&\leq& \frac{ \kappa_0}{ ( |x+b|-R)^{d+\alpha}}  \e^{-\lambda} 
\sum_{n=1}^\infty\frac{n^{d+\alpha}}{n!}
\left(\sum_{i=1}^n\int_{\mR^{nd}}1_{\sum_{j=1}^n z_j\in B_R(x+b)}\nu(\dif z_1)\cdots\dif z_i\cdots\nu(\dif z_n)\right)\\
&=& \frac{ \kappa_0}{ (|x+b|-R)^{d+\alpha}}  \e^{-\lambda}|B_R|\sum_{n=1}^\infty
\frac{n^{d+\alpha+1}\lambda^{n-1}}{n!} \\
&\leq& \frac{c_1 (d, \kappa_0) \kappa_0}{ (|x+b|-R)^{d+\alpha}}  \, |B_R|
\leq \frac{c_2 (d, \kappa_0, \alpha, R)}{(1+|x|)^{d+\alpha}}.
\end{eqnarray*}
where in the second equality we used the translation invariance property of the Lebesgue measure.
On the other hand, for any $x\in \R^d$, since $R>2$,
\begin{align*}
&\mP(\chi_2\in B_R(x))\geq\e^{-\lambda}\int_{\mR^{d}}1_{\{z_1\in B_R(x+b)\}}\nu(\dif z_1)\\
&\quad\geq \frac{ \kappa_0^{-1}\e^{-\lambda}|B_R(x+b)\cap B^c_1| }{(|x|+|b|+R)^{d+\alpha}}\\
&\quad\geq 
\frac{c_3 (d, \kappa_0, \alpha, R)}{(1+|x|)^{d+\alpha}}.
\end{align*}
Combining the above calculations, we get the desired estimate.
\end{proof}

Now we can give

\begin{proof}[Proof of Theorem \ref{Th1}]
Our proof is adapted from \cite{Wa}. Let $R>2$.
For the lower bound, by (i) of Lemma \ref{Le22}, we have 
$$
\delta:=\inf_{z\in B_R}p_1 (z)>0.
$$
Hence, by \eqref{ER5} and Lemma \ref{Le23},  
$$
p^\kappa_{0,1}(x)=\E \left[ p_1(x-\chi_2)\right]\geq \delta\mP(|x-\chi_2|\leq R) \geq \delta c^{-1}_1(1+|x|)^{-d-\alpha}.
$$
For the upper bound, by \eqref{ER5} again, we have
\begin{align}\label{ER6}
p^\kappa_{0,1}(x)\leq\mE \Big(p_1(x-\chi_2)1_{|x-\chi_2|\leq|x|/2}\Big)+\sup_{|z|>|x|/2}p_1(z).
\end{align}
By (ii) of Lemma \ref{Le22}, we can choose $N$-points $z_1,\cdots,z_N\in B_{|x|/2}$ such that 
$$
\mbox{$B_{|x|/2}\subset \cup_{j=1}^NB_\eps(z_j)$ and }\sum_{j=1}^N\sup_{z\in B_\eps(z_j)}p_1 (z)\leq c_4,
$$
where $c_4$ only depends on $\eps,\kappa_0,d,\alpha$.
Hence, by Lemma \ref{Le23}, we have
\begin{align*}
\mE \Big(p_1 (x-\chi_2)1_{|x-\chi_2|\leq|x|/2}\Big)
&\leq\sum_{j=1}^N\mE \Big(p_1 (x-\chi_2)1_{x-\chi_2\in B_\eps(z_j)}\Big)\\
&\leq\sum_{j=1}^N\sup_{z\in B_\eps(z_j)}p_1 (z)\mP(|x-\chi_2-z_j|\leq \eps)\\
&\leq c_1\sum_{j=1}^N\sup_{z\in B_\eps(z_j)}p_1 (z)(1+|x-z_j|)^{-d-\alpha}\\
&\leq c_1c_4(1+|x|/2)^{-d-\alpha},
\end{align*}
which together with  
 Lemma \ref{Le22}(ii)  yields 
\eqref{ER11}
 for $t=0$ and $s=1$. The proof is complete by \eqref{Sca}.
\end{proof}

Below for a function $f:\mR^d\to\mR$, $\theta\in(0,2)$ and $t>0$, we write
$$
\delta^{(\theta)}_f(t,x;z):=f(x+z)-f(x)-
\left(1_{\theta\in(1,2)}+1_{\theta=1}1_{|z|<t} \right)z\cdot\nabla f(x).
$$
The role of the variable $t$ in the definition of $\delta^{(\theta)}_f$ for $\theta=1$ is that we need to use the scaling 
property \eqref{Sca} for 
heat kernels
 $p^\kappa_{t, s}$. Note that under conditions \eqref{Ass1a} and \eqref{Ass1b}, 
$$
\sL^\kappa_t f(x)= \int_{\mR^d} \delta^{(\alpha)}_f (t, x; z) \frac{\kappa (t, x; z)}{|z|^{d+\alpha}} \dif z.
$$
By Theorem \ref{Th1}, the following lemma can be shown by the same convolution technique used in \cite[Lemma 2.3]{Ch-Zh}.

\bl\label{Le24}
Under \eqref{Ass1a}-\eqref{Ass1b}, for any $\theta\in(0,2)$,
there is a constant $c=c(\kappa_0,d,\alpha,\theta)>0$ such that
 for  every $s>t>0$ and $x,x',z\in\mR^d$,
 \begin{align}
|p^\kappa_{t,s}(x)-p^\kappa_{t,s} (x')|\leq c(((s-t)^{-1/\alpha}|x-x'|)\wedge 1)\Big(\varrho^0_\alpha(s-t,x)
+\varrho^0_\alpha(s-t,x')\Big),\label{ER222}
\end{align}
\begin{align}
|\nabla p^\kappa_{t,s}(x)|\leq c\varrho^0_{\alpha-1}(s-t,x),\label{ER70}
\end{align}
and 
\begin{align}
|\delta^{(\theta)}_{p^\kappa_{t,s} }(s-t,x;z)|\leq c\ell^{(\theta)}((s-t)^{-1/\alpha}z)
\left( \varrho^0_\alpha(s-t,x+z)+\varrho^0_\alpha(s-t,x)\right),\label{ER22}
\end{align}
\begin{align}
& |\delta^{(\theta)}_{p^\kappa_{t,s} }(s-t,x;z)-\delta^{(\theta)}_{p^\kappa_{t,s} }(s-t,x';z)|
 \leq c \left(((s-t)^{-1/\alpha}|x-x'|)\wedge 1\right) \ell^{(\theta)}((s-t)^{-1/\alpha}z)\no\\
&\qquad \times\left( \varrho^0_\alpha(s-t,x+ z)+\varrho^0_\alpha(s-t,x)+\varrho^0_\alpha(s-t,x'+ z)+\varrho^0_\alpha(s-t,x')\right),\label{ER33}
\end{align}
where $\ell^{(\theta)}(z):=1_{\theta\in(1,2)}(|z|^2\wedge|z|)+1_{\theta=1}(|z|^2\wedge 1)+1_{\theta\in(0,1)}(|z|\wedge 1)$.
\el

Using this lemma, it is easy to derive by definition (see \cite[Theorem 2.4]{Ch-Zh}).

\bl\label{Th24}
Under \eqref{Ass1a}-\eqref{Ass1b}, for any $\theta\in(0,2)$, there is a constant $c=c(\kappa_0,d,\alpha,\theta)>0$ such that
for all $s>t>0$ and $x,x'\in\mR^d$,
 \begin{align}
\int_{\mR^d}|\delta^{(\theta)}_{p^\kappa_{t,s} }(s-t,x;z)|\cdot|z|^{-d-\theta}\dif z&\leq c\varrho^0_{\alpha-\theta}(s-t,x),\label{ER6}
\end{align}
and
\begin{align}\label{ER7}
\begin{split}
&\int_{\mR^d}|\delta^{(\theta)}_{p^\kappa_{t,s} }(s-t,x;z)-\delta^{(\theta)}_{p^\kappa_{t,s} }(s-t,x';z)|\cdot |z|^{-d-\theta}\dif z\\
&\qquad\leq c(((s-t)^{-1/\alpha}|x-x'|)\wedge 1)\Big\{\varrho^0_{\alpha-\theta}(s-t,x)+\varrho^0_{\alpha-\theta}(s-t,x')\Big\}.
\end{split}
\end{align}
\el
Next we show the continuous dependence of $p^\kappa_{t,s}(x)$ with respect to $\kappa$.

\bl\label{Le26}
Let $\kappa$ and $\tilde\kappa$ be two kernels satisfying \eqref{Ass1a}-\eqref{Ass1b} with the same constant $\kappa_0$. 
Let $\alpha\in(0,2)$ and $\gamma\in[0,\alpha-1)$ when $\alpha\in(1,2)$, $\gamma\in[0,\alpha)$ if $\alpha\in(0,1]$. Assume that for some $K>0$,
\begin{align}\label{JG}
|\kappa(t,z)-\tilde\kappa(t,z)|\leq K\big(|z|^\gamma+1\big).
\end{align}
Then for any $\theta\in(0,2)$, there exists a constant $c=c(\kappa_0,d,\alpha,\theta)>0$ such that
\begin{align}
|\nabla^j p^{\kappa}_{t,s}(x)-\nabla^j p^{\tilde\kappa}_{t,s}(x)|\leq 
cK \varrho^0_{\alpha-j}(s-t,x),\ \ j=0,1,\label{ER308}
\end{align}
\begin{align}
\int_{\mR^d}|\delta^{(\theta)}_{p^\kappa_{t,s} }(s-t,x;z)-\delta^{(\theta)}_{p^\kappa_{t,s} }(s-t,x;z)|\cdot |z|^{-d-\theta}\dif z
\leq cK\varrho^0_{\alpha-\theta}(s-t,x).\label{ER77}
\end{align}
\el

\begin{proof}
By the scaling property \eqref{Sca}, it suffices to prove \eqref{ER308} and \eqref{ER77} for $s=1$ and $t=0$. 
The argument of proving \eqref{ER308} in \cite[Theorem 2.5]{Ch-Zh} strongly depends on the symmetry of $\kappa(z)=\kappa(-z)$.
Here we provide a different proof. Noticing that by \eqref{Ch} and \eqref{Den1},
$$
q_\lambda(x):=
p^{\lambda\kappa+(1-\lambda)\tilde\kappa}_{0,1}(x) 
=\int_{\mR^d}\e^{-\mathrm{i}x\cdot\xi}
\exp\left\{\int_{\mR^d}(\e^{\mathrm{i}\xi\cdot z}-1-\mathrm{i}\xi\cdot z^{(\alpha)})
\tfrac{\lambda\kappa(z)+(1-\lambda)\tilde\kappa(z)}{|z|^{d+\alpha}}\dif z\right\}\dif\xi,
$$
where $\kappa(z):=\int^1_0\kappa(r,z)\dif r$ 
and $\tilde\kappa(z):=\int^1_0\tilde\kappa(r,z)\dif r$.
We claim that 
\begin{equation}\label{e:2.26}
\p_\lambda q_\lambda(x)=(\sL^\kappa-\sL^{\tilde\kappa})q_\lambda(x).
\end{equation} 
In fact, by Fourier's transform it suffices to show that
$$
\p_\lambda \hat q_\lambda(\xi)=\widehat {\sL^\kappa q_\lambda}(\xi)-\widehat{\sL^{\tilde\kappa}q_\lambda}(\xi).
$$
Notice that 
$$
\widehat {\sL^\kappa f}(\xi)=\left(\int_{\mR^d}(\e^{\mathrm{i}\xi\cdot z}-1-\mathrm{i}\xi\cdot z^{(\alpha)})
\tfrac{\kappa(z)}{|z|^{d+\alpha}}\dif z\right)\hat f(\xi)
$$
and
$$
\hat q_\lambda(\xi)=\exp\left\{\int_{\mR^d}(\e^{\mathrm{i}\xi\cdot z}-1-\mathrm{i}\xi\cdot z^{(\alpha)})
\tfrac{\lambda\kappa(z)+(1-\lambda)\tilde\kappa(z)}{|z|^{d+\alpha}}\dif z\right\}.
$$
The desired claim now follows.
Now, by \eqref{e:2.26} and the assumption,
\begin{align}
|p^\kappa_{0,1}(x)-p^{\tilde\kappa}_{0,1}(x)|
&=\left|\int^1_0\p_\lambda q_\lambda(x)\dif \lambda\right|
=\left|\int^1_0\!\!\!\int_{\mR^d}\delta^{(\alpha)}_{q_\lambda}(1,x;z)
\frac{\kappa(z)-\tilde\kappa(z)}{|z|^{d+\alpha}}\dif z\dif \lambda\right|\no\\
&\leq K\int^1_0\!\!\!\int_{\mR^d}\big|\delta^{(\alpha)}_{q_\lambda}(1,x;z)\big| 
\frac{|z|^\gamma+1}{|z|^{d+\alpha}}\dif z\dif \lambda.\label{ER69}
\end{align}
If  $\alpha\in(1,2)$, since $\gamma\in[0,\alpha-1)$, by definition we have
$$
\int_{\mR^d}\big|\delta^{(\alpha)}_{q_\lambda}(1,x;z)\big|\cdot|z|^{\gamma-d-\alpha}\dif z
=\int_{\mR^d}\big|\delta^{(\alpha-\gamma)}_{q_\lambda}(1,x;z)\big|\cdot|z|^{\gamma-d-\alpha}\dif z.
$$
If $\alpha=1$, then
$$
\int_{\mR^d}\big|\delta^{(1)}_{q_\lambda}(1,x;z)\big|\cdot|z|^{\gamma-d-1}\dif z
\leq\int_{\mR^d}\big|\delta^{(1-\gamma)}_{q_\lambda}(1,x;z)\big|\cdot|z|^{\gamma-d-1}\dif z
+|\nabla q_\lambda(x)|\int_{|z|\leq 1}|z|^{\gamma-d}\dif z.
$$
If $\alpha\in(0,1)$, then
$$
\int_{\mR^d}\big|\delta^{(\alpha)}_{q_\lambda}(1,x;z)\big|\cdot|z|^{\gamma-d-\alpha}\dif z
=\int_{\mR^d}\big|\delta^{(\alpha-\gamma)}_{q_\lambda}(1,x;z)\big|\cdot|z|^{\gamma-d-\alpha}\dif z.
$$
Thus, by \eqref{ER69} and \eqref{ER6}, we obtain
\begin{align}\label{UY0}
|p^{\kappa (z)}_{0,1}(x)-p^{\tilde\kappa (z)}_{0,1}(x)|\leq cK(\varrho^0_{\alpha-\gamma}(1,x)+\varrho^0_{0}(1,x))=cK\varrho^0_{\alpha}(1,x),
\end{align}
which gives \eqref{ER308} for $j=0$. 
As for \eqref{ER308} with $j=1$ and \eqref{ER77}, it follows by the same argument as used in \cite[Theorem 2.5]{Ch-Zh}.
\end{proof}

\begin{remark}\label{R:2.7} \rm 
The assumption \eqref{JG} will be used to treat the more general kernel functions as in \eqref{Con1}. For instance, 
 $\kappa(t,x,z) :=\cos(x\cdot z)+2$ 
satisfies \eqref{Con1}, but is not H\"older continuous in $x$ uniformly with respect to $z$.
\end{remark}

\begin{remark}\label{R:2.7} \rm 
 The strong Markov process $X^{\kappa}_{t, s}$ of \eqref{e:2.3a} has infinitesimal
generator
$$
\sL^{\kappa}_t f(x) =  \int_{\R^d}
\left( f(x+z)-f(x)-\nabla f (x)\cdot z^{(\alpha)} \right) \frac{\kappa (t, z)}{|z|^{d+\alpha}} \dif z.
$$
The reason we use $z^{(\alpha)}$ instead of the more common $z 1_{\{|z| \leq 1\}}$
in the first order correction term is that 
this is the form 
for general $\alpha$-stable L\'evy processes.
Suppose
\begin{equation}\label{e:2.27}
\wt X^\kappa_{t,s}:=\int^s_t\!\!\!\int_{\mR^d}z \tilde N(\dif r,\dif z)+\int^s_t\!\!\!\int_{\{z\in \mR^d: |z|>1\}} z\, \frac{\kappa(r,z)}{|z|^{d+\alpha}}\dif z\dif r,
\end{equation} 
which has infinitesimal generator 
$$ 
\wt \sL^{\kappa}_t f(x) :=  \int_{\R^d}
\left( f(x+z)-f(x)-\nabla f (x)\cdot z 1_{\{|z|\leq 1\}} \right) \frac{\kappa (t, z)}{|z|^{d+\alpha}} \dif z.
$$
Clearly,
\begin{equation}\label{e:2.28}
\wt X^\kappa_{t,s} =  X^\kappa_{t,s} + \int_t^s b(r) \dif r
\quad {\mbox{and}} \quad 
\wt \sL^{\kappa}_t =\sL^{\kappa}_t + b(r) \cdot \nabla, 
\end{equation} 
where  
\begin{equation}\label{e:2.29}
b(r)=
\begin{cases}
- \int_{\{z\in \mR^d: |z|\leq 1\}} z \, \frac{\kappa (r, z)}{|z|^{d+\alpha}} \dif z
&\hbox{if } \alpha <1 , \\
0 & \hbox{if } \alpha =1 , \\
\int_{\{z\in \mR^d: |z|> 1\}} z \, \frac{\kappa (r, z)}{|z|^{d+\alpha}} \dif z
&\hbox{if } \alpha >1 .
\end{cases}
\end{equation}
 Denote by $\wt p^{\kappa}_{t, s}(x)$ 
the density function of $\wt X^\kappa_{t,s}$.
Then by \eqref{e:2.28},
$$
 \wt p^{\kappa}_{t, s}(x) =  p^{\kappa}_{t, s} \Big(x - \int_t^s b(r) \dif r \Big).
$$
Thus under conditions \eqref{Ass1a}-\eqref{Ass1b},
one can get two-sided estimates on $\wt p^{\kappa}_{t, s}(x)$ from that of 
$p^{\kappa}_{t, s}(x) $. 
In particular, it follows from Theorem \ref{Th1} that for every 
$\alpha \in [1, 2)$, 
there is a constant $ c_0>1$ depending only
on $\kappa_0$, $d$ and $\alpha$ so that for every $t<s$ and $x\in \mR^d$, 
$$ 
c_0^{-1} \frac{ s-t}{((s-t)^{1/\alpha} + |x|)^{d+\alpha}} 
\leq \wt p^\kappa_{t, s} (x) \leq c_0 \frac{ s-t}{((s-t)^{1/\alpha} + |x|)^{d+\alpha}} . 
$$
However, the above estimates in general fails for $\wt p^\kappa_{t, s} (x)$ when $\alpha 
\in (0, 1)$ as the drift $\int_t^s b(r) \dif r$ may not be controlled by $(s-t)^{1/\alpha}$ 
when $0<s-t\leq 1$.
\end{remark}  

\bigskip

We close this section by giving a proof to Theorem  \ref{Th01}.

\medskip

\begin{proof}[Proof of Theorem \ref{Th01}]
 First note that by \eqref{UR1} and Fubini's theorem,
$$
\sL^{1/2}f(x) =\frac{1}{2\Gamma(1/2)}\int_{\mR^d}(f(x-z)-f(x))\left(\int^\infty_0p_t(z)t^{-3/2}\dif t\right)\dif z,
$$
where $p_t(z)$ is the heat kernel of $\sL$ and so $P_tf(x)=\int_{\mR^d}f(z)p_t(x-z)\dif z$ .
   Note also that $\sL^{1/2}$ is 
  a L\'evy-type operator with L\'evy measure $\tilde \nu(\dif z):=\frac{1}{2\Gamma(1/2)}\int^\infty_0p_t(z)t^{-3/2}\dif t\dif z$.
By Theorem \ref{Th1}, one sees that
$$
\int^\infty_0p_t(z)t^{-3/2}\dif t\asymp\int^\infty_0t/(t^{1/\alpha}+|z|)^{d+\alpha}t^{-3/2}\dif t=|z|^{-d-\alpha/2}\int^\infty_0t^{-1/2}/(t^{1/\alpha}+1)^{d+\alpha}\dif t.
$$
Hence
$$
\tilde\nu(\dif z)=\tilde\kappa(z)\dif z/|z|^{d+\alpha/2} \
\hbox{ with $\tilde\kappa(z)\asymp 1$,}
$$
and by \cite[Corollary 4.4]{Zh}, for any $p>1$,
\begin{align}\label{EU1}
\|\sL^{1/2} f\|_p\asymp\|\Delta^{\alpha/4}f\|_p.
\end{align}
On the other hand, by \cite{St}, it is well known that for any $p>2d/(d+2\alpha)$,
\begin{align}\label{EU2}
\|\Delta^{\alpha/4}f\|_p\asymp\left\|\left(\int_{\mR^d}\frac{(f(\cdot+z)-f(\cdot))^2}{|z|^{d+\alpha}}\dif z\right)^{1/2}\right\|_p.
\end{align}
Moreover, by \eqref{Con3}, it is clear that 
\begin{align}\label{EU3}
\|\Gamma(f)^{1/2}\|_p\asymp\left\|\left(\int_{\mR^d}\frac{(f(\cdot+z)-f(\cdot))^2}{|z|^{d+\alpha}}\dif z\right)^{1/2}\right\|_p.
\end{align}
The desired estimate \eqref{Riesz} follows by combining \eqref{EU1}-\eqref{EU3}.
\end{proof}

\section{Proof of Theorem \ref{heat}}

In this section we consider the space and time dependent nonlocal operator
$\sL^{\kappa}_t$ defined by \eqref{Non}, with the kernel
function $\kappa (t,x, z)$ satisfying conditions \eqref{Con1}-\eqref{Con2},
and give a proof for Theorem \ref{heat}. 
 In order to emphasize the dependence on $x$, we also write
$$
\sL^{\kappa(x)}_t f(x)=\sL^{\kappa}_t f(x)=\int_{\mR^d} \delta^{(\alpha)}_f(t,x;z)\kappa(t,x,z)|z|^{-d-\alpha}\dif z.
$$
We use Levi's construction.  For fixed $y\in\mR^d$, let $\sL^{\kappa(y)}_t$ be the freezing operator
$$
\sL^{\kappa(y)}_t f(x)=\int_{\mR^d} \delta^{(\alpha)}_f(t,x;z)\kappa(t,y,z)|z|^{-d-\alpha}\dif z.
$$
Let $p^{(y)}_{t,s}(x):=p^{\kappa(y)}_{t,s}(x)$ be the heat kernel of operator $\sL^{\kappa(y)}_t$, i.e.,
\begin{align}
\p_t p^{(y)}_{t,s}(x)+\sL^{\kappa(y)}_t p^{(y)}_{t,s}(x)=0,\quad \lim_{t\uparrow s}p^{(y)}_{t,s}(x)=
\delta_{\{0\}} (x),  \label{ES2}
\end{align}
 where  $\delta_{\{ 0\}} (x)$ denotes the usual Dirac function.

Now, we want to seek the heat kernel $p^\kappa_{t,s} (x,y)$ of $\sL^{\kappa}_t$ with the following form:
\begin{align}
p^\kappa_{t,s} (x,y)=p^{(y)}_{t,s}(x-y)+\int^s_t\!\!\!\int_{\mR^d}p^{(z)}_{t,r}(x-z)q_{r,s}(z,y)\dif z\dif r.\label{ER65}
\end{align}
The classical Levi's method
suggests that $q_{t,s}(x,y)$ solves the following integral equation:
\begin{align}
q_{t,s}(x,y)=q^{(0)}_{t,s}(x,y)+\int^s_t\!\!\!\int_{\mR^d}q^{(0)}_{t,r}(x,z)q_{r,s}(z,y)\dif z\dif r,\label{EU2}
\end{align}
where
$$
q^{(0)}_{t,s}(x,y):=(\sL^{\kappa(x)}_{t}-\sL^{\kappa(y)}_{t})p^{(y)}_{t,s}(x-y).
$$
In fact, we formally have
\begin{align}
\p_tp^\kappa_{t,s} (x,y)&=\int^s_t\!\!\!\int_{\mR^d}\p_tp^{(z)}_{t,r}(x-z)q_{r,s}(z,y)\dif z\dif r-q_{t,s}(x,y)-\sL^{\kappa(y)}_{t}p^{(y)}_{t,s}(x-y)\no\\
&=-\int^s_t\!\!\!\int_{\mR^d}\sL^{\kappa(x)}_t p^{(z)}_{t,r}(x-z)q_{r,s}(z,y)\dif z\dif r-\sL^{\kappa(x)}_{t}p^{(y)}_{t,s}(x-y)\no\\
&=-\sL^{\kappa(x)}_{t}p^\kappa_{t,s} (\cdot,y)(x).\label{ER91}
\end{align}
The main point for us is to make the above calculations rigorous.
First of all, by \eqref{Con1}, we have for any $\eta\in(0,1)$,
$$
|\kappa(t,x,z)-\kappa(t,y,z)|\leq \kappa_0|x-y|^{\beta\eta}(1+|z|^{\beta'\eta}).
$$
Thus, in the following, without loss of generality, we may assume that
$$
\beta'\in(0,\alpha-1) \mbox{ if $\alpha\in(1,2)$ and $\beta'\in(0,\alpha)$ if $\alpha\in(0,1]$}.
$$
Noticing that by definition and \eqref{Con2},
\begin{align*}
q^{(0)}_{t,s}(x,y)=\int_{\mR^d}\delta^{(\alpha)}_{p^{(y)}_{t,s}}(s-t,x-y;z)(\kappa(t,x,z)-\kappa(t,y,z))|z|^{-d-\alpha}\dif z,
\end{align*}
as in \eqref{UY0}, we get by \eqref{Con1} and \eqref{ER6},
\begin{align}\label{JH1}
|q^{(0)}_{t,s}(x,y)|\leq c\varrho^\beta_0(s-t,x-y).
\end{align}
Hence,  by Picard's iteration, we can show (see \cite[Theorem 3.1]{Ch-Zh})
\bt\label{T3.4}
For $n\in\mN$, define $q^{(n)}_{t,s}(x,y)$ recursively by
\begin{align}
q^{(n)}_{t,s}(x,y):=\int^s_t\!\!\!\int_{\mR^d}q^{(0)}_{t,r}(x,z)q^{(n-1)}_{r,s}(z,y)\dif z\dif r.\label{EU22}
\end{align}
Under (\ref{Con1}) and (\ref{Con2}),  the series $q_{t,s}(x,y):=\sum_{n=0}^{\infty}q^{(n)}_{t,s}(x,y)$ is
 absolutely and locally uniformly convergent on $\mD^\infty_0$
and solves the integral equation (\ref{EU2}). Moreover,
 $q_{t,s}(x,y)$ is jointly continuous in $\mD^\infty_0$, and has the following estimates: For any $T>0$, 
there is a constant $c_1=c_1 (T,\kappa_0, d, \alpha, \beta)>0$ so that on $\mD^T_0$,
\begin{align}
|q_{t,s}(x,y)| \leq c_1(\varrho^\beta_0+\varrho^0_\beta)(s-t,x-y), \label{eq3}
\end{align}
and for any $\gamma\in(0,\beta)$,
there is a constant $c_2=c_2 (T, \kappa_0, d, \alpha, \beta, \gamma)>0$ so that on $\mD^T_0$,
\begin{align}\label{eq4}
\begin{split}
&|q_{t,s}(x,y)-q_{t,s}(x',y)| \leq c_2
\left( |x-x'|^{\beta-\gamma}\wedge 1\right)\\
 &\times\Big( (\varrho^0_\gamma+\varrho^\beta_{\gamma-\beta})(s-t,x-y)+(\varrho^0_\gamma+\varrho^\beta_{\gamma-\beta})(s-t,x'-y)\Big) .
\end{split}
\end{align}
\et

\medskip

Now using Lemmas \ref{Le24}-\ref{Le26} and as in \cite{Ch-Zh} (see also \cite{CHXZ}), 
we can make the calculations in \eqref{ER91} rigorous, and 
show that $p^\kappa_{t, s}(x, y)$ defined by \eqref{ER65} has the properties stated in Theorem \ref{heat} except for \eqref{Grad} and \eqref{GR2}.   
Note that in Theorem \ref{heat}, property (ii)   is implied 
 by property (iv).
 Below we  give a proof for  \eqref{Grad} and \eqref{GR2},
which are slight extensions of \cite[Lemma 4.2]{Ch-Zh}.

\begin{proof}[Proofs of \eqref{GR2} and \eqref{Grad}]
We only show \eqref{Grad} since \eqref{GR2} is similar by replacing $\Delta^{1/2}$ with $\nabla$ and using \eqref{ER70} and \eqref{ER308}. 
Let $\theta\in (0,(\alpha+\beta)\wedge 2)$. 
By \eqref{ER65}, we have
\begin{align*}
\Delta^{\theta/2} p^\kappa_{t,s} (x,y)&=\Delta^{\theta/2} p^{(y)}_{t,s}(x-y)+
\int^s_{(s+t)/2}\!\int_{\mR^d}\Delta^{\theta/2} p^{(z)}_{t,r}(x-z) q_{r,s}(z,y)\dif z\dif r\\
&\quad+\int^{(s+t)/2}_t\!\!\!\int_{\mR^d}\Delta^{\theta/2} p^{(z)}_{t,r}(x-z)(q_{r,s}(z,y)-q_{r,s}(x,y))\dif z\dif r\\
&\quad+\int^{(s+t)/2}_t\left(\int_{\mR^d}\Delta^{\theta/2} p^{(z)}_{t,r}(x-z)\dif z\right) q_{r,s}(x,y)\dif r\\
&=:J_1+J_2+J_3+J_4.
\end{align*}
For $J_1$, we have by \eqref{ER6} 
$$
I_1\lesssim\varrho^0_{\alpha-\theta}(s-t,x-y).
$$
For $J_2$, by \eqref{ER6}, \eqref{eq3} and \cite[Lemma 2.1]{Ch-Zh}, we have
\begin{align*}
J_2\lesssim \int^s_{(s+t)/2}\!\int_{\mR^d}\varrho^0_{\alpha-\theta}(r-t,x-z)(\varrho^0_\beta+\varrho^\beta_0)(s-r,z-y)\dif z\dif r\lesssim
\varrho^0_{\alpha-\theta}(s-t,x-y).
\end{align*}
For $J_3$, we get by \eqref{ER6}, \eqref{eq4} and \cite[Lemma 2.1]{Ch-Zh} that, for any $\gamma\in(0,(\alpha+\beta-\theta)\wedge\beta)$, 
\begin{align*}
J_3&\lesssim \int^{(s+t)/2}_t\!\!\!\int_{\mR^d}\varrho^{\beta-\gamma}_{\alpha-\theta}(r-t,x-z)(\varrho^0_\gamma+\varrho^\beta_{\gamma-\beta})(s-r,z-y)\dif z\dif r\\
&+\int^{(s+t)/2}_t\!\!\!\int_{\mR^d}\varrho^{\beta-\gamma}_{\alpha-\theta}(r-t,x-z)(\varrho^0_\gamma+\varrho^\beta_{\gamma-\beta})(s-r,x-y)\dif z\dif r\\
&\lesssim(\varrho^0_{\alpha+\beta-\theta}+\varrho^{\beta-\gamma}_{\alpha+\gamma-\theta}+\varrho^\beta_{\alpha-\theta})(s-t,x-y)\\
&+\int^{(s+t)/2}_t\!\!\!(r-t)^{(\beta-\gamma-\theta)/\alpha}(\varrho^0_\gamma+\varrho^\beta_{\gamma-\beta})(s-r,x-y)\dif r\\
&\lesssim\varrho^0_{\alpha-\theta}(s-t,x-y).
\end{align*}
For $J_4$, noticing that by \eqref{ER77},
\begin{align*}
\left|\int_{\mR^d}\Delta^{\theta/2} p^{(z)}_{t,r}(x-z)\dif z\right|
&=\left|\int_{\mR^d}(\Delta^{\theta/2} p^{(z)}_{t,r}-\Delta^{\theta/2} p^{(x)}_{t,r})(x-z)\dif z\right|\\
&\lesssim\int_{\mR^d}\varrho^\beta_{\alpha-\theta}(r-t,x-z)\dif z\lesssim (r-t)^{(\beta-\theta)/\alpha},
\end{align*}
we have by \eqref{eq3},
$$
I_4\lesssim\int^{(s+t)/2}_t(r-t)^{(\beta-\theta)/\alpha}(\varrho^0_\beta+\varrho^\beta_0)(s-r,x-y)\dif r\lesssim
\varrho^0_{\alpha-\theta}(s-t,x-y).
$$
Combining the above calculations, we get \eqref{Grad}.
\end{proof}

\medskip

\begin{proof}[Proof of Uniqueness]$\!\!$
We now show the uniqueness. 
Let $\wt p^\kappa_{t,s}(x,y)$ be another continuous function on $\mD^\infty_0$ satisfying 
\eqref{EQ}-\eqref{cz1}.
For $f\in C_c(\mR^d)$, a continuous function with compact support,
let $\wt P^\kappa_{t,s}f(x):= \int_{\mR^d}p^\kappa_{t,s}(x,y)f(y)\dif y$. 
By \eqref{EQ} and \eqref{cz1}, one sees that $\tilde u(t,x):=\wt P^\kappa_{t,s}f(x)$ solves the following equation
$$
\p_t \tilde u+\sL^\kappa_t \tilde u=0 \quad \hbox{with} \quad \lim_{t\uparrow s}\|\tilde u(t)-f\|_\infty=0.
$$
Note that by (i)-(ii) of Theorem \ref{heat}, $u(t, x):= P^\kappa_{t,s}f(x)$ has the same property as that for $\tilde u(t, x)$. 
Hence by the maximum principle (can be proved in a similar way as that for \cite[Theorem 6.1]{CHXZ}),  we have 
$\tilde u(t, x)=u(t, x)$; that is, 
$P^\kappa_{t,s}f(x)=\wt P^\kappa_{t,s}f(x)$. This implies $p^\kappa_{t,s}(x,y)=\wt p^\kappa_{t,s}(x,y)$.
\end{proof}

\end{document}